\documentclass[12pt,reqno]{article}

\usepackage[usenames]{color}
\usepackage{amssymb}
\usepackage{amsmath}
\usepackage{amsthm}
\usepackage{amsfonts}
\usepackage{amscd}
\usepackage{graphicx}

\usepackage[colorlinks=true,
linkcolor=webgreen,
filecolor=webbrown,
citecolor=webgreen]{hyperref}

\definecolor{webgreen}{rgb}{0,.5,0}
\definecolor{webbrown}{rgb}{.6,0,0}

\usepackage{color}

\usepackage{graphics}
\usepackage{latexsym}

\begin{document}

\theoremstyle{plain}
\newtheorem{theorem}{Theorem}
\newtheorem{corollary}[theorem]{Corollary}
\newtheorem{lemma}[theorem]{Lemma}
\newtheorem{proposition}[theorem]{Proposition}

\theoremstyle{definition}
\newtheorem{definition}[theorem]{Definition}
\newtheorem{example}[theorem]{Example}
\newtheorem{conjecture}[theorem]{Conjecture}

\theoremstyle{remark}
\newtheorem{remark}[theorem]{Remark}

\begin{center}
\vskip 1cm{\LARGE\bf Factors of Alternating Convolution of the Gessel Numbers
\vskip 1cm}
\large
Jovan Miki\'{c}\\
University of Banja Luka\\
Faculty of Technology\\
Bosnia and Herzegovina\\
\href{mailto:jovan.mikic@tf.unibl.org}{\tt jovan.mikic@tf.unibl.org} \\
\end{center}

\vskip .2in

\begin{abstract}
The Gessel number $P(n,r)$ is the number of the paths in plane with $(1,0)$ and $(0,1)$ steps from $(0,0)$ to $(n+r, n+r-1)$ that never touch any of the points from the set  $\{(x,x) \in \mathbb{Z}^2: x \geq r\}$. We show that there is a close relationship between the Gessel numbers $P(n,r)$ and the super Catalan numbers $S(n,r)$.
By using new sums, we prove that an alternating convolution of the Gessel numbers $P(n,r)$ is always divisible by $\frac{1}{2}S(n,r)$.

\end{abstract}

\noindent\emph{ \textbf{Keywords:}} Gessel Number, Super Catalan Number, Catalan Number, $M$ sum, Stanley's formula.

\noindent \textbf{2020} {\it \textbf{Mathematics Subject Classification}}:
05A10, 11B65.

\section{Introduction}\label{l:1}

Let $n$ be a non-negative integer, and let $r$ be a fixed positive integer. Let the number $P(n,r)$ denote $\frac{r}{2(n+r)}\binom{2n}{n}\binom{2r}{r}$. We shall call  $P(n,r)$ as the $n$th Gessel number of order $r$.

It is known that $P(n,r)$ is always an integer. They have an interesting combinatorial interpretation. The Gessel number $P(n,r)$  \cite[p.\ 191]{IG} counts all paths in the plane with unit horizontal and vertical steps from $(0, 0)$ to $(n+r, n+r-1)$ that never touch 
any of the points $(r,r)$, $(r+1, r+1)$, $\ldots$, . 

Recently \cite[Theorem 5, p.\ 2]{JM6}, it is shown that $P(n,r)$ is the number of all lattice paths in plane with $(1,0)$ and $(0,1)$ steps from $(0,0)$ to $(n+r, n+r-1)$ that never touch any of the points from the set $\{(x,x)\in \mathbb{Z}^2:1\leq x \leq n\}$; where $n$ and $r$ are positive integers.

Let $C_n=\frac{1}{n+1}\binom{2n}{n}$ denote the $n$th Catalan number, and let  $S(n,r)=\frac{\binom{2n}{n}\binom{2r}{r}}{\binom{n+r}{n}}$ denote the $n$th super Catalan number of order $r$. 
Obviously,  $P(n,1)=C_n$. Furthermore, it is readily verified that
\begin{equation}\label{eq:1}
P(n,r)=\binom{n+r-1}{n}\frac{1}{2}S(n,r)\text{.}
\end{equation}

It is known that  $S(n,r)$ is always an even integer except for the case $n=r=0$. See \cite[Introduction]{EAllen} and \cite[Eq.~(1), p.\ 1]{DC}. For only a few values of $r$, there exist combinatorial interpretations of $S(n,r)$. See, for example \cite{EAllen, DC, ChenWang,PipSch, Sch}. The problem of finding a combinatorial interpretation for super Catalan numbers of an arbitrary order $r$ is an intriguing open problem.

By using Eq.~(\ref{eq:1}),  it follows that  $P(n,r)$ is an integer.  Note that Gessel numbers $P(n,r)$ have a generalization \cite[Eq.~(1.10), p.\ 2]{VG2}.
Also, it is known   \cite[p.\ 191]{IG} that, for a fixed positive integer $r$, the smallest positive integer $K_r$ such that $\frac{K_r}{n+r}\binom{2n}{n}$ is an integer for every $n$ is $\frac{r}{2}\binom{2r}{r}$.

Let us consider the following sum:
\begin{equation}
\varphi(2n,m,r-1)=\sum_{k=0}^{2n}(-1)^k\binom{2n}{k}^m P(k,r)P(2n-k,r)\text{;}\label{eq:2}
\end{equation}
where $m$ is a positive integer.

For $r=1$, the sum in Eq.~(\ref{eq:2}) reduces to 
\begin{equation}
\varphi(2n,m,0)=\sum_{k=0}^{2n}(-1)^k\binom{2n}{k}^m C_k C_{2n-k}\text{.}\label{eq:3}
\end{equation}
Recently, by using new method,  it is shown \cite[Cor.\ 4, p.\ 2]{JM4} that the sum $\varphi(2n,m,0)$ is divisible by $\binom{2n}{n}$ for all non-negative integers $n$ and for all positive integers $m$. In particular, $\varphi(2n,1,0)=C_n\binom{2n}{n}$. See \cite[Th.\ 1, Eq.~(2)]{JM5}.  Interestingly, Gessel numbers appear \cite[Eq.~(68), p.\ 17]{JM4} in this proof.

By using the Eq.~(\ref{eq:1}), the sum $ \varphi(2n,m,r-1)$ can be rewritten, as follows:
\begin{equation}\label{eq:4}
\sum_{k=0}^{2n}(-1)^k\binom{2n}{k}^m \binom{k+r-1}{k}\binom{2n-k+r-1}{2n-k}\frac{1}{2}S(k,r)\frac{1}{2}S(2n-k,r)\text{.}
\end{equation}
 
Let $\varPsi(2n,m,r-1)$ denote the following sum 
\begin{equation}\label{eq:5}
\sum_{k=0}^{2n}(-1)^k\binom{2n}{k}^m S(k,r)S(2n-k,r)\text{.}
\end{equation}

Recently, by using new method, it is shown \cite[Th.\ 3, p.\ 3]{JM3} that the sum $\varPsi(2n,m,r-1)$ is divisible by $S(n,r)$ for all non-negative integers $n$ and $m$. In particular, $\varPsi(2n,1,r-1)=S(n,r)S(n+r,n)$. See \cite[Th.\ 1, Eq.~(1), p.\ 2]{JM3}.

Also,  it is known \cite[Th.\ 12]{JM2} that the sum $\sum_{k=0}^{2n}(-1)^k\binom{2n}{k}^m \binom{k+r-1}{k}\binom{2n-k+r-1}{2n-k}$ is divisible by $lcm(\binom{2n}{n},\binom{n+r-1}{n})$ for all non-negative integers $n$ and all positive integers $m$ and $r$.

Our main result is, as follows:
\begin{theorem}\label{t:1}
The sum $\varphi(2n,m,r-1)$ is always divisible by $\frac{1}{2}S(n,r)$ for all non-negative integers $n$ and for all positive integers $m$ and $r$.
\end{theorem}

For proving Theorem \ref{t:1}, we use a new class of binomial sums \cite[Eqns.~(27) and (28)]{JM2} that we call $M$ sums.
\begin{definition}\label{def:1}
Let $n$ and $a$ be  non-negative integers, and let $m$ be a positive integer. Let $S(n,m,a)=\sum_{k=0}^{n}\binom{n}{k}^m F(n,k,a)$, where $F(n,k,a)$ is an integer-valued function. Then  $M$ sums for the sum $S(n,m,a)$ are, as follows:
\begin{equation}\label{eq:6}
M_S(n,j,t;a)=\binom{n-j}{j}\sum_{v=0}^{n-2j}\binom{n-2j}{v}\binom{n}{j+v}^t F(n, j+v, a)\text{;}
\end{equation}
where $j$ and $t$ are non-negative integers such that $j\leq \lfloor \frac{n}{2} \rfloor$.
\end{definition}

Obviously, the following equation \cite[Eq.~(29)]{JM2},
\begin{equation}\label{eq:7}
S(n,m,a)=M_S(n,0,m-1;a)
\end{equation}
holds.

Let  $n$, $j$, $t$, and $a$ be same as in Definition \ref{def:1}.

It is known \cite[Th.\ 8]{JM2} that  $M$ sums satisfy the following recurrence relation:
\begin{equation}\label{eq:8}
M_S(n,j,t+1;a)=\binom{n}{j}\sum_{u=0}^{\lfloor\frac{n-2j}{2}\rfloor}\binom{n-j}{u}M_S(n,j+u,t;a)\text{.}
\end{equation}

Moreover,  we shall prove that $M$ sums satisfy another interesting recurrence relation:
\begin{theorem}\label{t:2}
Let  $n$ and $a$ be non-negative integers, and let $m$ be a positive integer. Let $R(n,m,a)$ denote 
\begin{equation}
\sum_{k=0}^{n}\binom{n}{k}^m G(n,k,a)\text{;}\notag
\end{equation}
where $G(n,k, a)$ is an integer-valued function. Let $Q(n,m,a)$ denote
\begin{equation}
\sum_{k=0}^{n}\binom{n}{k}^m H(n,k,a)\text{;}\notag
\end{equation}
where $H(n,k,a)=\binom{a+k}{a}\binom{a+n-k}{a}G(n,k,a)$.
Then the following recurrence relation is true:
\begin{equation}
M_Q(n,j,0;a)=\binom{a+j}{a}\sum_{l=0}^{a}\binom{n-j+l}{l}\binom{n-j}{a-l}M_R(n,j+a-l,0;a)\text{.}\label{eq:9}
\end{equation}
\end{theorem}

Note that, by using substitution $k=v+j$, the Eq.~(\ref{eq:6}) becomes, as follows:
\begin{equation}
M_S(n,j,t;a)=\binom{n-j}{j}\sum_{k=j}^{n-j}\binom{n-2j}{k-j}\binom{n}{k}^t F(n, k, a)\text{.}\label{eq:10}
\end{equation}\label{eq:15}
From now on, we use the Eq.~(\ref{eq:10}) instead of the Eq.~(\ref{eq:6}).

\section{Background}\label{l:2}
In $ 1998$, Calkin  proved that the alternating binomial sum $S_1(2n, m)=\sum_{k=0}^{2n}(-1)^k\binom{2n}{k}^m $ is divisible by $\binom{2n}{n} $ for all non-negative integers $n$ and all positive integers $m$.   In 2007, Guo, Jouhet, and Zeng proved, among other things, two generalizations of Calkin's result \cite[Thm.\ 1.2, Thm.\ 1.3, p.\ 2]{VG1}. In $2018$, Calkin's result \cite[Thm.\ 1]{NC} is proved by using $D$ sums \cite[Section 8]{JM1}. Note that there is a close relationship between $D$ sums and $M$ sums \cite[Eqns.~(9) and (19)]{JM2}.

Recently, Calkin's result \cite[Thm.\ 1]{NC} is proved by using $M$ sums \cite[Section 5]{JM2}. In particular, it is known \cite[Eqns.~(22) and (25)]{JM2} that: 

\begin{align*}
M_{S_1}(2n,j,0)&=
\begin{cases}0, & \text{if } 0 \leq j <n;\\
(-1)^n, & \text{if } j=n.
\end{cases}\\
M_{S_1}(2n,j,1)&=(-1)^n\binom{2n}{n}\binom{n}{j}\text{.}
\end{align*}

By using $M$ sums, it is shown \cite[Section 7]{JM2}  that the sum $S_2(2n,m,a)\\=\sum_{k=0}^{2n}(-1)^k\binom{2n}{k}^m\binom{a+k}{k}\binom{a+2n-k}{2n-k}$ is divisible by $lcm(\binom{a+n}{a}, \binom{2n}{n})$ for all non-negative integers $n$ and $a$; and for all positive integers $m$.
In particular, it is known  \cite[Eqns.~(52) and (57)]{JM2} that:
\begin{align*}
M_{S_2}(2n,j,0;a)&=(-1)^{n-j}\binom{a+n}{a}\binom{a+j}{j}\binom{a}{n-j},\\
M_{S_2}(2n,j,1;a)&=(-1)^{n-j}\binom{2n}{n}\sum_{u=0}^{n-j}(-1)^u\binom{n}{j+u}\binom{j+u}{u}\binom{a+j+u}{j+u}\binom{a+n}{2n-j-u}\text{.}
\end{align*}

By using $D$ sums, it is shown \cite[Thm.\ 1]{JM4} that the sum $S_3(2n,m)\\ =\sum_{k=0}^{2n}(-1)^k\binom{2n}{k}^m\binom{2k}{k}\binom{2(2n-k)}{2n-k}$ is divisible by $\binom{2n}{n}$ for all non-negative integers $n$, and for all positive integers $m$.  The same result also can be proved by using $M$ sums. It can be shown that 
\begin{equation}\label{eq:11}
M_{S_3}(2n,j,0)=(-1)^j\binom{2n}{n}\binom{2j}{j}\binom{2(n-j)}{n-j}\text{.}
\end{equation}
Note that Eq.~(\ref{eq:11}) is equivalent with the \cite[Eq.~(12)]{JM4}.

Recently, by using $D$ sums, it is shown  \cite[Th.\ 3, p.\ 3]{JM3}  that $\varPsi(2n,m,l-1)$  is divisible by $S(n,l)$ for all non-negative integers $n$ and $m$. The same result also can be proved by using $M$ sums. Let $l$ be a positive integer, and let $n$ and $j$ be non-negative integers such that $j\leq n$. It is readily verified \cite[Eq.~(91)]{JM2} that
\begin{equation}\label{eq:12}
M_{\varPsi}(2n,j,0;l-1)=(-1)^j\frac{\binom{2l}{l}\binom{2n}{n}\binom{2j}{j}\binom{2(n+l-j)}{n+l-j}\binom{2n-j}{n}}{\binom{n+l}{n}\binom{2n+l-j}{n}}\text{.}
\end{equation}

Furthermore, by using \cite[Eq.~(103)]{JM3} and \cite[Eq.~(33)]{JM2}, it can be shown that
\begin{equation}\label{eq:13}
M_{\varPsi}(2n,j,1;l-1)=(-1)^jS(n,l)\binom{n}{j}\sum_{v=0}^{n-j}(-1)^vS(n+l-j-v,n)\binom{2(j+v)}{j+v}\binom{n-j}{v}\text{.}
\end{equation}

The rest of the paper is structured as follows. In Section \ref{l:3}, we give a proof of Theorem \ref{t:2} by using Stanley's formula.
In Section \ref{l:4}, we give a proof of Theorem \ref{t:1}. Our proof of Theorem \ref{t:1} consists from two parts. In the first part, we prove that Theorem \ref{t:1} is true for $m=1$. In the second part, we prove that Theorem \ref{t:1} is true for all postive integers $m$ such that $m\geq 2$.

\section{A Proof of Theorem \ref{t:2}}\label{l:3}

We use three known binomial formulas.

Let $a$, $b$, and $c$ be non-negative integers such that $a \geq b \geq c$. The first formula is:
\begin{equation}
\binom{a}{b}\binom{b}{c}=\binom{a}{c}\binom{a-c}{b-c}\text{.}\label{eq:14}
\end{equation}

Let $a$, $b$, $m$, and $n$ be non-negative integers. The second formula is Stanley's formula:
\begin{equation}
\sum_{k=0}^{min(m,n)}\binom{a}{m-k}\binom{b}{n-k}\binom{a+b+k}{k}=\binom{a+n}{m}\binom{b+m}{n}\text{.}\label{eq:15}
\end{equation}

The third formula is a symmetry of binomial coefficients.

\begin{proof}

By setting  $S:=Q$ and $t:=0$  in the Eq.~(\ref{eq:10}), it follows that
\begin{align}
M_Q(n,j,0;a)&=\binom{n-j}{j}\sum_{k=j}^{n-j}\binom{n-2j}{k-j}\binom{n}{k}^{0} H(n,k,a)\text{,}\notag\\
&=\sum_{k=j}^{n-j}\binom{n-j}{j}\binom{n-2j}{k-j}\binom{a+k}{a}\binom{a+n-k}{a}G(n,k,a)\label{eq:16}\text{.}
\end{align}

By using  symmetry of binomial coefficients and the Eq.~(\ref{eq:14}), it follows that
\begin{equation}
\binom{n-j}{j}\binom{n-2j}{k-j}=\binom{n-j}{k-j}\binom{n-k}{n-k-j}\text{.}\label{eq:17}
\end{equation}

By using the Eq.~(\ref{eq:17}), Eq.~(\ref{eq:16}) becomes, as follows:
\begin{align}
M_Q(n,j,0;a)&=\sum_{k=j}^{n-j}\binom{n-j}{k-j}\binom{n-k}{n-k-j}\binom{a+k}{a}\binom{a+n-k}{a}G(n,k,a)\notag\text{,}\\
&=\sum_{k=j}^{n-j}\binom{n-j}{k-j}\bigl{(}\binom{a+n-k}{a}\binom{n-k}{n-k-j}\bigr{)}\binom{a+k}{a}G(n,k,a)\text{.}\label{eq:18}
\end{align}

By using  symmetry of binomial coefficients and the Eq.~(\ref{eq:18}), it follows that
\begin{equation}
\binom{a+n-k}{a}\binom{n-k}{n-k-j}=\binom{a+n-k}{a+j}\binom{a+j}{j}\text{.}\label{eq:19}
\end{equation}

By using the Eq.~(\ref{eq:19}), the Eq.~(\ref{eq:18}) becomes, as follows:
\begin{equation}
M_Q(n,j,0;a)=\binom{a+j}{j}\sum_{k=j}^{n-j}\binom{n-j}{k-j}\bigl{(}\binom{a+n-k}{a+j}\binom{a+k}{a}\bigr{)}G(n,k,a)\text{.}\label{eq:20}
\end{equation}

By setting $m:=a+j$, $n:=a$, $a:=n-k$, $b:=k-j$, and $k:=l$ in Stanley's formula \ref{eq:15}, we obtain that
\begin{equation}\label{eq:21}
\binom{a+n-k}{a+j}\binom{a+k}{a}=\sum_{l=0}^{a}\binom{n-k}{a+j-l}\binom{k-j}{a-l}\binom{n-j+l}{l}\text{.}
\end{equation}

By using  Eqns.~(\ref{eq:20}) and (\ref{eq:21}), it follows that  $M_Q(n,j,0;a)$ is equal to 
\begin{equation}
\binom{a+j}{j}\sum_{k=j}^{n-j}\binom{n-j}{k-j}\bigl{(}\sum_{l=0}^{a}\binom{n-k}{a+j-l}\binom{k-j}{a-l}\binom{n-j+l}{l}\bigr{)}G(n,k,a)\text{.}\label{eq:22}
\end{equation}

Note that, by the Eq.~(\ref{eq:21}), we obtain the following inequality:
\begin{align}
n-k\geq a+j-l\text{, or }\notag\\
k\leq n-j-a+l\text{.}\label{neq:1}
\end{align}
Similarly,  by the Eq.~(\ref{eq:21}), we obtain another inequality
\begin{align}
k-j\geq a-l\text{, or}\notag\\
k\geq a+j-l\text{.}\label{neq:2}
\end{align}

By changing the order of sumation in the Eq.~(\ref{eq:22}) and by using Inequalites (\ref{neq:1}) and (\ref{neq:2}), it follows that $M_Q(n,j,0;a)$ is equal to 
\begin{equation}
\binom{a+j}{j}\sum_{l=0}^{a}\binom{n-j+l}{l}\sum_{k=a+j-l}^{n-(a+j-l)}\binom{n-j}{k-j}\binom{n-k}{a+j-l}\binom{k-j}{a-l}G(n,k,a)\text{.}\label{eq:25}
\end{equation}

By using symmetry of binomial coefficients and the Eq.~(\ref{eq:14}), it follows that
\begin{equation}
\binom{n-j}{k-j}\binom{n-k}{a+j-l}=\binom{n-j}{a+j-l}\binom{n-a-2j+l}{k-j}\text{.}\label{eq:26}
\end{equation}

By using Eqns.~(\ref{eq:22}) and (\ref{eq:26}), it follows that $M_Q(n,j,0;a)$ is equal to 

\begin{equation}
\binom{a+j}{j}\sum_{l=0}^{a}\binom{n-j+l}{l}\binom{n-j}{a+j-l}\sum_{k=a+j-l}^{n-(a+j-l)}\binom{n-a-2j+l}{k-j}\binom{k-j}{a-l}G(n,k,a)\text{.}\label{eq:27}
\end{equation}

By using symmetry of binomial coefficients and the Eq.~(\ref{eq:14}), it follows that
\begin{equation}
\binom{n-a-2j+l}{k-j}\binom{k-j}{a-l}=\binom{n-a-2j+l}{a-l}\binom{n-2a-2j+2l}{k-j-a+l}\text{.}\label{eq:28}
\end{equation}

By using  Eqns.~(\ref{eq:27}) and (\ref{eq:28}), it follows that $M_Q(n,j,0;a)$ is equal to 

\begin{equation}
\binom{a+j}{j}\sum_{l=0}^{a}\binom{n-j+l}{l}\binom{n-j}{a+j-l}\binom{n-a-2j+l}{a-l}\sum_{k=a+j-l}^{n-(a+j-l)}\binom{n-2a-2j+2l}{k-j-a+l}G(n,k,a)\text{.}\label{eq:29}
\end{equation}

Again, by using symmetry of binomial coefficients and the Eq.~(\ref{eq:14}), it follows that
\begin{equation}
\binom{n-j}{a+j-l}\binom{n-a-2j+l}{a-l}=\binom{n-j}{a-l}\binom{n-(j+a-l)}{j+a-l}\text{.}\label{eq:30}
\end{equation}

By using  Eqns.~(\ref{eq:29}) and (\ref{eq:30}), it follows that $M_Q(n,j,0;a)$ is equal to 

\begin{equation}
\binom{a+j}{j}\sum_{l=0}^{a}\binom{n-j+l}{l}\binom{n-j}{a-l}\bigl{(}\binom{n-(j+a-l)}{j+a-l}\sum_{k=a+j-l}^{n-(a+j-l)}\binom{n-2(a+j-l)}{k-(j+a-l)}G(n,k,a)\bigr{)}\text{.}\label{eq:31}
\end{equation}

Note that , by setting $S:=R$, $F:=G$, $j:=j+a-l$, and $t:=0$ in the Eq.~(\ref{eq:10}), it follows that:
\begin{equation}
M_R(n,j+a-l,0;a)=\binom{n-(j+a-l)}{j+a-l}\sum_{k=a+j-l}^{n-(a+j-l)}\binom{n-2(a+j-l)}{k-(j+a-l)}G(n,k,a)\text{.}\label{eq:32}
\end{equation}

Hence, by using  Eqns.~(\ref{eq:31}) and (\ref{eq:32}), it follows that 
\begin{equation}
M_Q(n,j,0;a)=\binom{a+j}{j}\sum_{l=0}^{a}\binom{n-j+l}{l}\binom{n-j}{a-l}M_R(n,j+a-l,0;a)\text{.}\notag
\end{equation}

This completes the proof of Theorem \ref{t:2}.

\end{proof}

\section {A Proof of Theorem \ref{t:1}}\label{l:4}
Let  the function $\varphi(2n,m,r-1)$ be defined as in Eq.~(\ref{eq:4}).

Let us consider the following sum
\begin{equation}
\phi(2n,m,r-1)=\frac{1}{4}\varPsi(2n,m,r-1)\text{.}\label{eq:33}
\end{equation}
 By Eq.~(\ref{eq:5}), the Eq.~(\ref{eq:33}) becomes
\begin{equation}\label{eq:34}
\phi(2n,m,r-1)=\sum_{k=0}^{2n}(-1)^k\binom{2n}{k}^m( \frac{1}{2}S(k,r))(\frac{1}{2}S(2n-k,r))\text{.}
\end{equation}

Note that, since $r$ is a positive integer, both numbers $ \frac{1}{2}S(k,r)$ and $\frac{1}{2}S(2n-k,r)$ are integers. Therefore, the sum $\phi(2n,m,r-1)$ is a sum from Definition \ref{def:1}. 

Now we can apply Theorem \ref{t:2}. By setting $Q:=\varphi$, $n:=2n$, $a:=r-1$, $R:=\phi$ and  $G(2n,k,r-1):=\frac{1}{2}S(k,r)\frac{1}{2}S(2n-k,r)$ in Theorem \ref{t:2}, by the Eq.~(\ref{eq:9}), it follows that $M_{\varphi}(2n,j,0;r-1)$ is equal to

\begin{equation}\label{eq:35}
\binom{j+r-1}{j}\sum_{l=0}^{r-1}\binom{2n-j+l}{l}\binom{2n-j}{r-1-l}M_{\phi}(2n,j+r-1-l,0;r-1)\text{.}
\end{equation}

By  Eqns.~(\ref{eq:6}) and (\ref{eq:33}), it follows that
\begin{equation}\label{eq:36}
M_{\phi}(2n,j,t;r-1)=\frac{1}{4}M_{\varPsi}(2n,j,t;r-1)\text{.}
\end{equation}

By setting $t:=0$ and $l:=r$ in the Eq.~(\ref{eq:12}) and by using the Eq.~(\ref{eq:36}) and the definition of super Catalan numbers, we obtain that
\begin{equation}\label{eq:37}
M_{\phi}(2n,j,0;r-1)=(-1)^j\frac{1}{2}S(n, r)\frac{\binom{2j}{j}\binom{2(n+r-j)}{n+r-j}\binom{2n-j}{n}}{2\binom{2n+r-j}{n}}\text{.}
\end{equation}

By using Eq.~(\ref{eq:37}), it follows that the sum $M_{\phi}(2n,j+r-1-l,0;r-1)$ is equal to
\begin{equation}
(-1)^{j+r-1-l}\frac{1}{2}S(n,r)\frac{\binom{2(j+r-1-l)}{j+r-1-l}\binom{2(n-j+l+1)}{n-j+l+1}\binom{2n-j-r+l+1}{n}}{2\binom{2n-j+l+1}{n}}\text{.}
\end{equation}\label{eq:38}

By using a symmetry of binomial coefficients and the Eq.~(\ref{eq:14}), it follows that
\begin{equation}
\binom{2n-j}{r-1-l}\binom{2n-j-r+l+1}{n}=\binom{2n-j}{n}\binom{n-j}{j+r-l-1}\text{.}\label{eq:39}
\end{equation}

By using Eqns.~(38) and (39), it follows that the sum $\binom{2n-j}{r-1-l}M_{\phi}(2n,j+r-1-l,0;r-1)$ is equal to

\begin{equation}\label{eq:40}
(-1)^{j+r-1-l}\frac{1}{2}S(n,r)\binom{2n-j}{n}\frac{\binom{2(j+r-1-l)}{j+r-1-l}\binom{2(n-j+l+1)}{n-j+l+1}\binom{n-j}{j+r-l-1}}{2\binom{2n-j+l+1}{n}}\text{.}
\end{equation}

By using  Eqns.~(\ref{eq:35}) and (\ref{eq:40}), it follows that the sum $M_{\varphi}(2n,j,0;r-1)$ is equal to

\begin{equation}\label{eq:41}
(-1)^{j+r-1}\binom{j+r-1}{j}\frac{1}{2}S(n,r)\binom{2n-j}{n}\sum_{l=0}^{r-1}(-1)^l\binom{2n-j+l}{l}\frac{\binom{2(j+r-1-l)}{j+r-1-l}\binom{2(n-j+l+1)}{n-j+l+1}\binom{n-j}{j+r-l-1}}{2\binom{2n-j+l+1}{n}}\text{.}
\end{equation}

By setting $j:=0$ in the Eq.~(\ref{eq:41}), we obtain that
\begin{align}
M_{\varphi}(2n,0,0;r-1)&=(-1)^{r-1}\frac{1}{2}S(n,r)\binom{2n}{n}\sum_{l=0}^{r-1}(-1)^l\binom{2n+l}{l}\frac{\binom{2(r-1-l)}{r-1-l}\binom{2(n+l+1)}{n+l+1}\binom{n}{r-l-1}}{2\binom{2n+l+1}{n}}\notag\\
&=(-1)^{r-1}\frac{1}{2}S(n,r)\sum_{l=0}^{r-1}(-1)^l\binom{2n+l}{l}\binom{2(r-1-l)}{r-1-l}\binom{n}{r-l-1}\frac{\binom{2n}{n}\binom{2(n+l+1)}{n+l+1}}{2\binom{2n+l+1}{n}}\notag\\
&=(-1)^{r-1}\frac{1}{2}S(n,r)\sum_{l=0}^{r-1}(-1)^l\binom{2n+l}{l}\binom{2(r-1-l)}{r-1-l}\binom{n}{r-l-1}\bigl{(}\frac{1}{2}S(n,n+l+1)\bigr{)}\text{.}\label{eq:42}
\end{align}

By using the Eq.~(\ref{eq:41}), it follows that the sum $M_{\varphi}(2n,0,0;r-1)$ is divisible by $\frac{1}{2}S(n,r)$.  By setting $n:=2n$, $m:=1$, $S:=\varphi$, and $r:=r-1$ in the Eq.~(\ref{eq:7}), we obtain that
\begin{equation}\label{eq:43}
\varphi(2n,1,r-1)=M_{\varphi}(2n,0,0;r-1)\text{.}
\end{equation}

By using Eqns.~(\ref{eq:42}) and (\ref{eq:43}), it follows that the sum $\varphi(2n,1,r-1)$ is divisible by $\frac{1}{2}S(n,r)$. This completes the proof of Theorem \ref{t:1} for the case $m=1$.

Let us calculate the sum $M_{\varphi}(2n,j,1;r-1)$, where $n$ and $j$ are non-negative integers such that $j \leq n$.

By setting $n:=2n$, $S:=\varphi$, $t:=0$,  and $a:=r-1$ in the Eq.~(\ref{eq:8}), we obtain that
\begin{equation}\label{eq:44}
M_{\varphi}(2n,j,1;r-1)=\sum_{u=0}^{n-j}\binom{2n}{j}\binom{2n-j}{u}M_{\varphi}(2n,j+u,0;r-1)\text{.}
\end{equation}

By using a symmetry of binomial coefficients and the Eq.~(\ref{eq:14}), it follows that
\begin{equation}\label{eq:45}
\binom{2n}{j}\binom{2n-j}{u}=\binom{2n}{2n-j-u}\binom{j+u}{u}\text{.}
\end{equation}

By using Eqns.~(\ref{eq:44}) and (\ref{eq:45}), it follows that
\begin{equation}\label{eq:46}
M_{\varphi}(2n,j,1;r-1)=\sum_{u=0}^{n-j}\binom{j+u}{u}\bigl{(}\binom{2n}{2n-j-u}M_{\varphi}(2n,j+u,0;r-1)\bigr{)}\text{.}
\end{equation}

By using the Eq.~(\ref{eq:41}), it follows that the sum $ \binom{2n}{2n-j-u}M_{\varphi}(2n,j+u,0;r-1)$ equals 
\begin{equation}\label{eq:47}
\begin{split}
\binom{2n}{2n-j-u}\binom{2n-j-u}{n}(-1)^{j+u+r-1}\binom{j+u+r-1}{j+u}\frac{1}{2}S(n,r)\cdot\\\cdot\sum_{l=0}^{r-1}(-1)^l\binom{2n-j-u+l}{l}\frac{\binom{2(j+u+r-1-l)}{j+u+r-1-l}\binom{2(n-j-u+l+1)}{n-j-u+l+1}\binom{n-j-u}{j+u+r-l-1}}{2\binom{2n-j-u+l+1}{n}}\text{.}
\end{split}
\end{equation}

By using a symmetry of binomial coefficients and the Eq.~(\ref{eq:14}), it follows that
\begin{equation}\label{eq:48}
\binom{2n}{2n-j-u}\binom{2n-j-u}{n}=\binom{2n}{n}\binom{n}{j+u}\text{.}
\end{equation}

By using  Eqns.~(\ref{eq:47}) and (\ref{eq:48}), it follows that the sum  $ \binom{2n}{2n-j-u}M_{\varphi}(2n,j+u,0;r-1)$ equals

\begin{equation}
\begin{split}
\frac{1}{2}S(n,r)\binom{n}{j+u}(-1)^{j+u+r-1}\binom{j+u+r-1}{j+u}\sum_{l=0}^{r-1}(-1)^l\binom{2n-j-u+l}{l}\cdot\\\cdot\binom{2(j+u+r-1-l)}{j+u+r-1-l}\binom{n-j-u}{j+u+r-l-1}\bigl{(}\frac{1}{2}S(n,n-j-u+l+1)\bigr{)}\text{.}\label{eq:49}
\end{split}
\end{equation}

Hence, by using the Eq.~(\ref{eq:49}), we obtain that
\begin{equation}\label{eq:50}
\binom{2n}{2n-j-u}M_{\varphi}(2n,j+u,0;r-1)=\frac{1}{2}S(n,r)\cdot c(n,j+u,r-1)\text{;}
\end{equation}
where $c(n,j+u,r-1)$ is always an integer.

By using the Eq.~(\ref{eq:50}), the Eq.~(\ref{eq:46}) becomes, as follows:
\begin{equation}\label{eq:51}
M_{\varphi}(2n,j,1;r-1)=\frac{1}{2}S(n,r)\sum_{u=0}^{n-j}\binom{j+u}{u}c(n,j+u,r-1)\text{.}
\end{equation}

By Eq.~(\ref{eq:51}), it follows that the sum $M_{\varphi}(2n,j,1;r-1)$ is divisible by $\frac{1}{2}S(n,r) $ for all non-negative integers $n$ and $j$ such that $j\leq n$; and for all positive integers $r$.  By using the Eq.~(\ref{eq:8}) and the induction principle, it can be shown that
the sum $M_{\varphi}(2n,j,t;r-1)$ is divisible by $\frac{1}{2}S(n,r) $ for all non-negative integers $n$ and $j$ such that $j\leq n$; and for all positive integers $r$ and $t$. 

By setting $S=\varphi$, $n:=2n$, $m:=t+1$, and $a:=r-1$ in the Eq.~(\ref{eq:7}), it follows that

\begin{equation}\label{eq:51}
\varphi(2n,t+1,r-1)=M_{\varphi}(2n,0,t;r-1)\text{.}
\end{equation}

Since $t \geq 1$, it follows that $t+1\geq 2$. By Eq.~(\ref{eq:51}), it follows that the sum $\varphi(2n,m,r-1)$ is always divisible by $\frac{1}{2}S(n,r) $ for  all non-negative integers $n$; and for all positive integers $m$ and $r$ such that $m \geq 2$. See \cite[Section 4]{JM2}.
This completes the proof of Theorem \ref{t:1}

\begin{remark} \label{rem:1}
For $r=1$, by using Theorem \ref{t:1} and the fact $\frac{1}{2}S(n,1)=C_n$, it follows that the sum $\varphi(2n,m,0)$ is divisible by $C_n$. Therefore, for $r=1$, the result of Theorem \ref{t:1} is weaker than the result  that sum $\varphi(2n,m,0)$ is divisible by $\binom{2n}{n}$ \cite[Cor.\ 4, p.\ 2]{JM4}. However, for $n=3$, $r=2$, and $m=1$, the sum  $\varphi(2n,m,r-1)$ is  neither divisible by $\binom{2n}{n}$ nor by $S(n,r)$.
\end{remark}

\section{Conclusion}\label{l:5}
Theorem \ref{t:1} and Eq.~(\ref{eq:1}) show that there is a close relationship between Gessel's numbers $P(n,r)$ and super Catalan numbers $S(n,r)$. Gessel's numbers have, at least, two combinatorial interpretations \cite{IG, JM6}, while super Catalan numbers do not have a combinatorial interpretation in general. We think that some combinatorial properties of Gessel's numbers can be used for obtaining a combinatorial interpretation for super Catalan numbers in general.

\section*{Acknowledgments}
I want to thank professor Tomislav Do\v{s}li\'{c} for inviting me on the third Croatian Combinatorial Days in Zagreb.

\end{document}